\documentclass[11pt,a4paper]{article}

\usepackage{comment}
\usepackage{graphicx}
\usepackage{amsmath,bm}
\usepackage{amssymb}
\usepackage{amsfonts}
\usepackage{amsthm}
\usepackage{xcolor}
\usepackage[T1]{fontenc}
\usepackage{lmodern}

\numberwithin{equation}{section}

\usepackage{paralist}
\numberwithin{equation}{section}

\usepackage{mathtools, mathabx, verbatim}

\def\R{\mathbb{R}}

\def\Torus{\mathbb{T}^3}

\def\neweq#1{\begin{equation}\label{#1}}
\def\endeq{\end{equation}}

\newtheorem{thm}{Theorem}[section]
\newtheorem{lem}[thm]{Lemma}

\newtheorem{prop}[thm]{Proposition}

\newtheorem{defn}[thm]{Definition}
\theoremstyle{definition}
\newtheorem{rem}[thm]{Remark}
\theoremstyle{remark}

\newcommand{\ds}{\displaystyle}

\newcommand{\abs}[1]{\left\vert#1\right\vert}

\newcommand{\N}{\mathbb{N}}

\DeclareMathOperator{\dive}{div}

\DeclareMathOperator*{\esssup}{ess\,sup}
\DeclareMathOperator{\supp}{supp}

{\left\{\begin{array}{@{}l@{}}}{\end{array}\right.}

\makeatletter 
\def\@makefnmark{} 
\makeatother 

\title{\bf Velocity-vorticity geometric constraints for the energy conservation of 3D
  ideal incompressible fluids}
\author{Luigi C. Berselli and Rossano Sannipoli
  \\
  Dipartimento di Matematica - Università di Pisa, Italy
  \\
  email: luigi.carlo.berselli@unipi.it, rossano.sannipoli@dm.unipi.it}

\date{\today}
\begin{document}
\maketitle
\begin{abstract}
  In this paper we consider the 3D Euler equations and we first prove a criterion for
  energy conservation for weak solutions with velocity satisfying additional assumptions
  in fractional Sobolev spaces with respect to the space variables, balanced by proper
  integrability with respect to time. Next, we apply the criterion to study the energy
  conservation of solution of the Beltrami type, carefully applying properties of products
  in (fractional and possibly negative) Sobolev spaces and employing a suitable bootstrap
  argument.
  \\
  \textbf{MCS: }{Primary 35Q31; Secondary 76B03.}
\end{abstract}

\section{Introduction}
We consider the homogeneous incompressible 3D Euler equations
\begin{equation}\label{eq:Euler}
\begin{cases}
  \partial_t u + (u\cdot \nabla)\,u+\nabla p = 0 \qquad \qquad& (t,x) \in (0,T)\times
  \Torus, 
  \\
  \dive u = 0 & (t,x) \in (0,T)\times \Torus,
  \\
    u(0,x)=u_0(x) & x \in \Torus,
\end{cases}
\end{equation}
where $\Torus := \mathbb R^3 \backslash \mathbb Z^3$, and
$u: (0,T)\times \Torus\to \mathbb R^3$ and $p: (0,T)\times \Torus \to \mathbb R$ represent
respectively the velocity vector field and the kinematic pressure of an ideal fluid.  It
is well known that for smooth solution to~\eqref{eq:Euler} (which are known to exists only
locally in time) the kinetic energy
\begin{equation*}
    E(t) = \frac{1}{2}\int_{\Torus}\abs{u(t,x)}^2\,dx=\frac{1}{2}\|u(t)\|_2^2,
\end{equation*}
is constant. Let $u$ and $p$ smooth enough to perform the following calculations: we
rewrite the convective term of ~\eqref{eq:Euler}$_1$ as follows
$$
(u\cdot\nabla)\,u=\dive (u\otimes u).
$$
%
Multiplying~\eqref{eq:Euler}$_1$ by the solution itself, and integrating over the
domain, we get, that $\|u(t)\|_2^2=\|u_0\|_2^2$ since
  \begin{equation*}
    \int_{\Torus}\dive (u\otimes u):u\,dx=-    \int_{\Torus} (u\otimes u):\nabla u\,dx=-
    \int_{\Torus} u\cdot \nabla \frac{|u|^2}{2}\,dx=0.
  \end{equation*}
  We report this very basic calculation since we will use it several times and also since
  we will show how it changes with a curl-formulation of the convective term.

Since physical experiments show that taking the limit as the viscosity vanishes, the
energy dissipation seems not to vanish (cf.~Frisch~\cite{Fri1995}), it has been a subject
of many studies to understand if this conservation remains valid supposing certain
(limited) regularity on the solution of the Euler equation. In 1954, Lars
Onsager~\cite{Ons1949} conjectured that if $u$ is sufficiently regular in space, say
$u\in L^\infty(0,T;C^\theta)$, with $\theta >\frac{1}{3}$, then the kinetic energy is
preserved; on the other hand for $\theta <\frac{1}{3}$ a dissipation phenomenon could be
possible, even in absence of viscosity. The positive part of this conjecture was solved 40
years later by Constantin, E, Titi (see~\cite{CET1994}), where they proved a slightly more
general result in Besov spaces (which implies the H\"older case). See also
Eyink~\cite{Eyi1994}. To be more precise, in \cite{CET1994} it is proved the conservation
of energy if $u\in L^3(0,T;B^\theta_{3,\infty})$, for $\theta > \frac{1}{3}$. Sharpest
results were proved later on by Duchon and Robert~\cite{DR2000} and Cheskidov \textit{et
  al.}~\cite{CCFS2008}. Results in scales of classical H\"older functions are proved
in~\cite{Ber2023c}, while the boundary value problem is analyzed in Bardos and
Titi~\cite{BT2018}. In the last fifteen years --starting from the celebrated result by De
Lellis and  Sz{\'e}kelyhidi~\cite{DLS2009}-- also the negative part of the Onsager conjecture
has been addressed, with an endpoint in the work by Isett~\cite{Ise2018} and Buckmaster
\textit{et al.}~\cite{BDLSV2019}. Nevertheless there is still a strong activity to
determine the minimal space-time assumptions which are sufficient for the energy
conservation, and some recent results are those in~\cite{BG2024,WWWY2023}.

Taking inspiration also from the work by De Rosa~\cite{Der2020} and Liu, Wang, and
Ye~\cite{LWY2023}, we consider here criteria in scales of fractional Sobolev spaces,
instead of Besov or H\"older spaces. This will allow us also to obtain sharp results which
reach the critical exponents, see the discussion in Lemma~\ref{lem:lemlimsup}.

The first result we prove concerns the conservation of energy in the fractional Sobolev
setting. We restrict to the Hilbertian case $W^{s,2}(\Torus)=H^s(\Torus)$, but similar
results in scales of fractional Sobolev spaces $W^{s,p}(\Torus)$ can be obtained along the
same lines.

\begin{thm}\label{thm:energyconservationfracsob}
  Let $u\in L^\frac{5}{2s}(0,T, H^s(\Torus))$, with $\frac{5}{6}\le s< \frac{5}{2}$, be a
  weak solution to the Euler equation~\eqref{eq:Euler}. Then, the kinetic energy is
  conserved, that is
    \begin{equation*}
        \|u(t)\|_{L^2(\Torus)}= \|u_{0}\|_{L^2(\Torus)} \qquad \text{for a.e. } t \in [0,T].
    \end{equation*}
  \end{thm}
  We observe that the condition $L^3(0,T;H^{5/6}(\Torus))$ is exactly the same condition
  proved in Cheskidov, Friedlander, and Shvydkoy~\cite{CFS2010} for the Navier-Stokes
  equations (even if a more technical setting of the problem with boundaries), see also
  Beirao and Yiang~\cite[Prop.~4.5]{BY2019}, again in the viscous case.  In this paper we
  identify the same as a sufficient condition also for the problem without viscosity. In
  this respect note also that the extension to the Euler equations of results known for
  the Navier-Stokes equations is one of the results proved in~\cite{BG2024}. Note also
  that, similar to other observations in Nguyen, Nguyen, and Tang~\cite{NgNgT2019},
  Wang~\textit{et al.}~\cite{WWWY2023}, the criteria involve the ``critical spaces'' and
  not slightly smaller spaces, as when considering Besov or H\"older spaces,
  cf.~\cite{Ber2023c,CET1994}, where the smooth functions are not dense with respect to
  the norm of the space itself.

\medskip

As an application of the criteria in Theorem~\ref{thm:energyconservationfracsob} (which
are obviously valid also for periodic Leray-Hopf weak solutions to the Navier-Stokes
equations), we analyze the energy conservation of a family of solutions with a particular
geometric meaning: the Beltrami (also known as Trkal) flows. Beltrami solutions are well
known in fluid dynamics as they provide a family of stationary solutions to the Euler
equations~\eqref{eq:Euler}. These are such that the curl of the velocity field, denoted by
$\omega:=\nabla\times u$, is proportional to the field itself, i.e.
\begin{equation}\label{eq:Beltrami}
    \omega(x,t) = \lambda(x,t)u(x,t),
\end{equation}
where $\lambda(\cdot,\cdot)$ is a suitable scalar function of the space and/or time
variables.

Note that these flows, despite being in some cases very simple and smooth (note that for
instance potential flows are Beltrami flows with $\lambda\equiv0$) they are genuinely 3D
flows, since in 2D the (scalar) vorticity is orthogonal to the plane of motion.

In addition, we observe that, by using the so-called \textit{Lamb vector}
$\omega\times u$, it is possible to write the alternative \textit{rotational formulation}
of the convective term
\begin{equation*}
  (u\cdot\nabla)\, u=\omega\times u+\frac{1}{2}\nabla|u|^2.
\end{equation*}
This implies that in the case of Beltrami flows the convective term is equal to a gradient
(Bernoulli pressure) which can be included in the pressure: Beltrami flows (if smooth)
satisfy linear (non-local) evolution equations, since the quadratic term becomes simply a
gradient and the nonlinearity disappears. This means that the flow is laminar, but there
are two caveat: a) the numerical simulation of the pressure and especially that of the
Bernoulli pressure is particularly critical: if not using \textit{pressure-robust}
numerical methods, the result at very high Reynolds numbers could be affected by large
oscillations (see Gauger, Linke, and Schroeder~\cite{GLS2019}); b) more important from the
theoretical point of view is the fact that $(\omega\times u)\cdot u$ formally vanishes;
for non-smooth functions the fact $(\omega\times u)\cdot u=0$ (at most almost everywhere)
does not directly imply that
\begin{equation*}
  \int_0^t\int_{\Torus}(\omega\times u+\frac{1}{2}\nabla|u|^2)\cdot u\,dx d\tau=0,
\end{equation*}
since the integration could be not justified. One sufficient condition could be that of
showing that the above integral exists: then it will vanish, but unfortunately this is not
the case for weak solutions. At least with respect to the space variables,
$u \in L^2(\Torus)$ and the vorticity field is a distribution in $H^{-1}(\Torus)$ and so
the term $(\omega\times u)\cdot u$ could be not defined.
%

We start making some observations on the regularity which follows from the geometric
constraint~\eqref{eq:Beltrami}. If $\lambda(x,t)\equiv \lambda \in \mathbb R$ (such a
condition corresponds to the circularly polarized plane waves used in electromagnetism),
then $u$ is smooth and the conservation of energy is an obvious consequence.  This follows
by a standard bootstrap argument using the Biot-Savart formula: in fact using that
$-\Delta u=\text{curl }\omega$, from $u\in L^\infty(0,T;L^2(\Torus))$ we can infer by
elliptic regularity in the space variables that then
$\omega=\lambda u\in L^\infty(0,T;L^2(\Torus))$, which implies
$u\in L^\infty(0,T;H^1(\Torus))$. Iterating we get $u\in L^\infty(0,T;H^3)$, which is a
class of classical solutions. This implies that a continuation argument for smooth
solutions is valid, provided that the initial datum is smooth.

The second observation comes from a simple computation in the case in which
$\lambda(x,t)=\lambda(t)\in L^{p}(0,T)$, for some $p\ge 1$. Observe that in this case
$\nabla\cdot\omega=\lambda(t)\nabla\cdot u=0$, so the divergence-free constraint is
satisfied without any further assumption on $\lambda(t)$.  Then~\eqref{eq:Beltrami}
implies that $\omega \in L^p(0,T;L^2(\Torus))$, and consequently we have more regularity on $u$,
indeed $u\in L^p(0,T;H^1(\Torus))$. Iterating this procedure we get that if
$\lambda(t)\in L^p(0,T)$, for some $p\ge 1$, and $u\in L^\infty(0,T;L^2(\Torus))$ then 
\begin{equation*}
    \omega\in L^\frac{p}{3}(0,T;H^2(\Torus))\xhookrightarrow{}L^\frac{p}{3}(0,T;L^\infty(\Torus)).
\end{equation*}
Hence, if $p\geq3$,  this is the Beale-Kato-Majda~\cite{BKM1984} criterion for
continuation of smooth solutions (which conserve the energy). Hence, a first elementary
results is the following.
\begin{prop}
  Let $u$ be a weak solution to the Euler equations, which is a Beltrami solution with
  $\lambda\in L^p(0,T)$, with $p\geq3$. Let $u_0\in H^3(\Torus)$, with
  $\nabla\cdot u_0=0$, then $u$ is the unique a classical solution of~\eqref{eq:Euler} in
  $[0,T]$ and conserves the energy.
\end{prop}


If $\lambda$ depends also on the space variables a compatibility condition to preserve the
divergence condition is that $\nabla\lambda\cdot u=0$. This has some consequences on the
effective velocity fields to be considered, especially if the solutions are classical, see
Beltrami~\cite{Bel1873} and Trkal~\cite{Trk1919}.  For recent results on possible
existence of non-trivial Beltrami fields, see Enciso and Peralta~\cite{EP2016} and
Abe~\cite{Abe2022}. In our case we suppose to have a weak solution, which is also a
Beltrami field and work directly on it.

\medskip

The second Theorem we prove is about the conservation of energy when $u$ is a Beltrami
field.
\begin{thm}\label{thm:Beltramifieldthm}
  Let be a weak solution to the Euler
  equation~\eqref{eq:Euler}, such that it is a Beltrami field, i.e. \eqref{eq:Beltrami} is
  satisfied.   If $\lambda \in L^\beta(0,T,H^\tau(\Torus))$, with $\beta>\frac{5}{2\tau-1}$, for
    $\frac{1}{2}<\tau\leq\frac{3}{2}$, or if $\lambda\in L^{5/2}(0,T;H^\tau(\Torus))$, with 
    $\tau>\frac{3}{2}$, then the kinetic energy is conserved.
\end{thm}

This result derives from Theorem~\ref{thm:energyconservationfracsob} after a proper (even
iterated) use of some precise results about the continuity of the multiplication operator
in (negative) Sobolev spaces, see the results summarized in the next section.
\bigskip

\textbf{Plan of the paper:} In Section~\ref{sec:preliminaries} we give some basic
definition and introduce functional spaces we will work with, recalling different results
that will be used later on in the paper. In Section~\ref{sec:mainresults} we will give the
proofs of the two main results, namely
Theorems~\ref{thm:energyconservationfracsob}-\ref{thm:Beltramifieldthm}.
%

%
\section{Preliminaries}\label{sec:preliminaries}
\subsection{Functional Spaces and weak solutions}
In this paper we will use the classical periodic Lebesgue spaces
$(L^p(\mathbb{T}^3), \|\cdot\|_p)$ and Sobolev spaces
$(W^{k,p}(\mathbb{T}^3), \|\cdot\|_{W^{k,p}})$ of natural order $k\in\N$, all considered
with zero mean value. When $p=2$ we also use the notation
$H^{k}=W^{k,2}(\mathbb{T}^3)$. Here we will denote by $(\cdot,\cdot)$ and $\|\cdot\|$ the
scalar product and the norm in $L^2(\Torus)$ respectively. In addition we do not
distinguish norm of scalar or vector valued functions.

A central role in this paper will be played by the
fractional Sobolev spaces that we define in the following
(see~\cite{BH2021}).
\begin{defn}[Fractional Sobolev spaces]
  Let $s\in \mathbb R$ and $1\le p\le\infty$. We define the
  Sobolev-Slobodeckij spaces as follows
  \begin{itemize}
  \item Let $s\in (0,1)$, then we will say that
    $u\in W^{s,p}(\Torus) $ if
    \begin{equation*}
      \|u\|_{W^{s,p}}:= \|u\|_{p}+ [u]_{W^{s,p}} <\infty, 
    \end{equation*}
    where
    \begin{equation*}
      [u]_{W^{s,p}}=
      \begin{cases}
        \bigg(\ds\int_{\Torus}\int_{\Torus}
        \frac{|u(x)-u(y)|^p}{|x-y|^{n+s p}}\,dxdy\bigg)^\frac{1}{p} &
        p\in [1,+\infty)
        \\\\
        \ds\esssup\limits_{\substack{x,y\in \Torus
            \\
            x\neq y}}\frac{|u(x)-u(y)|}{|x-y|^{s}} & p=\infty;
      \end{cases}
    \end{equation*}
  \item Let $s=\theta + k$, with $\theta \in (0,1)$ and
    $k\in \mathbb N_0$. Then, we will say that
    $u\in W^{s,p}(\Torus)$ if
    \begin{equation*}
      \|u\|_{W^{s,p}}:= \|u\|_{W^{k,p}}+
      \sum_{|r|=k}[\partial_r u]_{W^{\theta, p} }<\infty; 
    \end{equation*} 
  \item If $s<0$, and $p'$ is the Sobolev conjugate exponent of $p$,
    then
    \begin{equation*}
      W^{s,p}(\Torus)=(W^{-s,p'}(\Torus))^*,
    \end{equation*}
    where $*$ denotes the topological dual space.
  \end{itemize}
%
\end{defn}
\begin{rem}
  Note that for $p=2$ one can use also an equivalent semi-norm
  $[u]_{W^{\alpha,2}} =\|(-\Delta)^{\alpha/2}u\|_{L^{2}}$.

  In the case of the whole space, but also for bounded domains with a proper definition of
  the restriction fractional spaces can be defined even by means of Bessel potentials, and
  in this case the space is denoted in literature by $H^{s,p}(\Omega)$. For our purposes,
  we do not give this definition since for $s\in \mathbb R$ and $p=2$, the two definitions
  coincide, see Triebel~\cite{Tri1978} . That is why we will denote $W^{s,2}(\Omega)$ as
  $H^s(\Omega)$.
\end{rem}
Here we state some propositions that will be useful when considering the product of two
fractional Sobolev functions $u \in H^{s_1}(\Torus)$, $\lambda \in H^{s_2}(\Torus)$
(see~\cite[Thm 6.1, 7.3, 8.1, 8.2]{BH2021} for the whole space case and the results in the
periodic setting follow along the same lines. Here results are rephrased in the simpler
case $p_i=p=2$ ).

The first proposition regards the case of non-negative exponents.
\begin{prop}\label{prop:prodnonnegative}
     Let $s, s_i \in \mathbb R$ be parameters such that for $i=1,2$
     \begin{enumerate}
    \item $s\ge 0$;
        \item $s_i \ge s$;
        \item $s_1+s_2-s > \frac32$.
    \end{enumerate}
    If $u\in H^{s_1}(\Torus)$ and $\lambda \in H^{s_2}(\Torus)$, then $\lambda u \in H^{s}(\Torus)$ and the map of pointwise multiplication
    \begin{equation*}
        H^{s_1}(\Torus)\times H^{s_2}(\Torus)\to H^{s}(\Torus),
    \end{equation*}
    is continuous and bilinear. Moreover if $s\in \mathbb N_0$, the strictness of inequalities $(2)$ and $(3)$ can be interchanged.
\end{prop}
In the case of non-negative exponents we have the following proposition
\begin{prop}\label{prop:prodnegative}
     Let $s, s_i \in \mathbb R$ be parameters such that for $i=1,2$
     \begin{enumerate}
    \item $s_i\ge s$;
        \item $\min\{s_1,s_2\}<0$;
        \item $s_1+s_2\ge 0$;
        \item $s_1+s_2-s> \frac32$.
    \end{enumerate}
    If $u\in H^{s_1}(\Torus)$ and $\lambda \in H^{s_2}(\Torus)$, then $\lambda u \in H^{s}(\Torus)$ and the map of pointwise multiplication
    \begin{equation*}
        H^{s_1}(\Torus)\times H^{s_2}(\Torus)\to H^{s}(\Torus),
    \end{equation*}
    is continuous and bilinear. Moreover, if we substitute the condition $(2)$ by the
    condition
    \begin{enumerate}
      \setcounter{enumi}{4}
        \item $\min\{s_1,s_2\}\ge0$ and $s<0$,
    \end{enumerate}
    and the inequality $(3)$ is strict, the same result holds.
\end{prop}
Additionally, if $s\notin \mathbb N_0$, proposition~\ref{prop:prodnonnegative} is valid
even in open and bounded set with Lipschitz boundary.

To complete these definitions we will say that the Bochner measurable
function $u\in L^q(0,T,W^{s,p}(\Omega)) $ if the
following norm
\begin{equation*}
  \|u\|_{L^q(W^{s,p})}:= 
    \begin{cases}
      \bigg(\int_0^T\|u(t)\|_{W^{s,p}}^q\,dt\bigg)^\frac{1}{q}\quad\text{ if }q<\infty,
      \\
      \esssup\limits_{t\in[0,T]} \|u(t)\|_{W^{s,p}}\quad\text{ if }q=\infty,
    \end{cases}
\end{equation*}
is finite.\\
We want to introduce now the notion of a weak solution to the Euler
equations, it is necessary to introduce the spaces $H$ and $V$, which
are respectively the closure in $L^2(\Torus)$ and
$W^{1,2}(\mathbb{T}^3)$ of the smooth, periodic, divergence-free and
with zero mean-value vector fields. The space of test functions will
be
\begin{equation*}
  \mathcal{D}_T = \{\varphi \in C_{0}^{\infty}([0,T[; C^\infty(\mathbb{T}^3)) :
  \nabla\cdot\varphi =0\}. 
\end{equation*}
Now we are able to define the notion of weak solutions for the Euler equations which we
will consider.
\begin{defn}
  Let $v_0\in H$. A measurable function $v:(0,T)\times \mathbb{T}^3\to \Torus$ is called a
  weak solution to the Euler equation if $v \in L^\infty(0,T,H)$ is such that
  \begin{equation*}
    \int_0^{\infty} [(v,\partial_t\varphi)+((v\otimes v),\nabla \varphi)]\,dt =
    -(v(0),\varphi(0)) \qquad  \forall\, \varphi \in \mathcal{D}_T. 
    \end{equation*}
\end{defn}
\subsection{Mollification}
A fundamental tool in the sequel will be that of mollification and we recall the most
relevant properties, stated for periodic functions. Let us consider a centrally symmetric
function $\rho \in C_0^\infty(\R^3)$, such that $\rho\ge 0$, $\supp{\rho}\subset B_1(0)$
and $\|\rho\|_{L^1(\R^3)}=1$. Let $\varepsilon\in (0,1]$, we define the family of
Friederichs mollifiers as follows
$\rho_\varepsilon(x)=\varepsilon^{-3}\rho(\varepsilon^{-1}x)$. Then for every
$f\in L^1_{\text{loc}}(\Torus)$ we can define the well-posed mollification of $f$, that is
\begin{equation*}
  f_\varepsilon(x) = \int_{\Torus}\rho_\varepsilon(x-y)f(y)\,dy =\int_{\Torus}\rho_\varepsilon(y)f(x-y)\,dy,
  \end{equation*}
  which is nothing else the convolution of $\rho_\varepsilon$ and $f$.
  \\
  Since for small $\varepsilon>0$,
  $\supp{\rho_\varepsilon}\subset B_\varepsilon(0)\subset ]-\pi,\pi[^3$, then we can say
  that
\begin{equation*}
    f_\varepsilon(x)=\int_{\Torus}\rho_\varepsilon(y)f(x-y)\,dy.
\end{equation*}
We note that if $f\in L^1(\Torus)$, then
$f\in L^1_{\text{loc}}(\Torus)$ and $f_\varepsilon$ is
$2\pi$-periodic along the $x_i$-axis, for $i=1,2,3$.

Apart classical result on mollification in Lebesgue, H\"older and, Sobolev spaces, most of
the results can be extended to fractional Sobolev spaces. An useful Lemma, contained
in~\cite{Der2020}, can be summarized as follows.
\begin{lem}
  Let $f,g: \Torus\to \Torus$ are such that $f\in W^{\alpha,p}(\Torus)$ and
  $g \in W^{\beta,q}(\Torus)$, for some $0<\alpha, \beta<1$ and $p,q \ge 1$ such that
  $\frac{1}{p}+\frac{1}{q}=\frac{1}{m}$. Then. for every $1\le m <+\infty$, there exists a
  constant $C=C(m)$, such that
    \begin{equation}\label{gradestimate}
        \|\nabla f_\varepsilon\|_{p} \le C \varepsilon^{\alpha-1}[f]_{W^{\alpha,p}},
      \end{equation}
    \begin{equation}\label{diffestimate}
        \|(f\otimes g)_\varepsilon-f_\varepsilon\otimes g_\varepsilon\|_{m}\le C
        \varepsilon^{\alpha+\beta}[f]_{W^{\alpha,p}}[g]_{W^{\beta,q}}. 
      \end{equation}
      Moreover, in the case $\alpha=\beta=0$ it also follows that
      \begin{equation}
        \label{eq:sharpestimate}
        \limsup_{\varepsilon\to0}\|(f\otimes g)_\varepsilon-f_\varepsilon\otimes g_\varepsilon\|_{m}=0.
      \end{equation}
\end{lem}
We end this section with a Lemma whose proof comes from the estimate~\eqref{gradestimate}
and the fact that $C^\infty(\Torus)$ is dense in $W^{\alpha,p}(\Torus)$, as it
is done in~\cite[Lemma 2.1]{NgNgT2019} for the case $\alpha=0$
\begin{lem}\label{lem:lemlimsup}
  Let $f: \Torus \to \Torus$ be a function in
  $L^q(0,T;W^{\alpha,p}(\Torus))$, for some $\alpha\in]0,1[$. 
  Then, for every  $1\le p,q <\infty$, we have
    \begin{equation*}
        \limsup_{\varepsilon \to 0^+}\varepsilon^{1-\alpha}\|\nabla f_\varepsilon\|_{L^q(L^p)}=0.
    \end{equation*}
    Moreover if $f \in L^q(0,T;L^{p}(\Torus))$, then
    \begin{equation*}
        \limsup_{\varepsilon \to 0^+}\|f_\varepsilon-f\|_{L^q(L^p)}=0.
    \end{equation*}    
  \end{lem}
\begin{proof}
  The proof of this result is based on the observation that, if $f\in
  W^{\alpha,p}(\Torus)$, then $f_{\varepsilon}$ is infinitely smooth but higher
  norms are not uniform bounded in $\varepsilon>0$. In particular,
  from~\eqref{gradestimate}  one can deduce immediately the
  boundedness 
  \begin{equation*}
    \sup_{\varepsilon>0}\varepsilon^{1-\alpha}\|\nabla
    f_{\varepsilon}\|_{L^{p}}\leq C[f]_{W^{\alpha,p}}<+\infty. 
  \end{equation*}
  To prove the Lemma it is enough to observe that for each $\lambda>0$ we can find
  $g\in C^\infty(\Torus)$ such that $\|f-g\|_{W^{\alpha,p}}<\lambda$. Then, applying again
  the estimate~\eqref{gradestimate} to $f-g$ one gets
  \begin{equation*}
    \varepsilon^{1-\alpha}\|\nabla f_{\varepsilon}\|_{p}\leq
    \varepsilon^{1-\alpha}\|\nabla g_{\varepsilon}\|_{p}+ C
    [f-g]_{W^{\alpha,p}}\leq     \varepsilon^{1-\alpha}\|\nabla g_{\varepsilon}\|_{p}+ C
    \lambda. 
  \end{equation*}
  Since $\lambda$ can be chosen arbitrarily small and since
  $\lim_{\varepsilon\to0}\varepsilon^{1-\alpha}\|\nabla g_{\varepsilon}\|_{p}=0$ (being
  $g$ smooth and fixed) we get the proof. The need for the density of smooth function
  excludes the case $p=\infty$ and, more generally excludes from this type of results
  Nikol'ski\u{\i} and H\"older spaces. Then, the extension to $f$ in the Bochner space
  $L^{q}(0,T;W^{\alpha,p}(\Torus))$ is simply obtained by raising to the $q$-th power the
  above estimate and integrating over $(0,T)$.
\end{proof}
%
%
\section{Main results}\label{sec:mainresults}
In this section we give the proof of the main results of the paper. We start with the
proof of a criterion about conservation of energy for velocities in the fractional Sobolev
spaces.
\begin{proof}[Proof of Theorem~\ref{thm:energyconservationfracsob}]
  By following a very standard procedure to deal with non-smooth functions we use as test
  function in the definition of weak solutions
  $\rho_{\varepsilon}*(\rho_{\varepsilon}*u)$.  To be precise the argument will also need
  another smoothing in time which is nevertheless standard to justify, see~\cite{BT2018}.
  By using the identity
  \begin{equation*}
    \int_0^t \int_{\Torus}(u_\varepsilon\otimes u_\varepsilon):\nabla u_\varepsilon\,dx
    d\tau=0,
  \end{equation*}
  being $u_{\varepsilon}$ smooth and divergence-free, we get the equality
  \begin{equation*}
    \frac{1}{2}    \|u_{\varepsilon}(t)\|^2-\frac{1}{2}\|u_{0,\varepsilon}\|^2=\int_0^t
    \int_{\Torus}((u\otimes u)_\varepsilon-u_\varepsilon\otimes u_\varepsilon):\nabla
    u_\varepsilon\,dx d\tau.
  \end{equation*}
  %
%
%
%
%
%
%
%
%
%
%
%
%
%
%
%
%
%
%
%
%
%
  If we manage to estimate the integrand of the previous inequality in such a way the
  right-hand side goes to zero, we have finished since the $L^2$-convergence of
  $u_{\varepsilon}(t)$ to $u(t)$ as $\varepsilon\to0$ holds almost everywhere in
  $t\in(0,T)$, by the properties of smoothing. The proof will be slightly different
  depending on the values of the exponent ``$s\in \R$'' in the extra-assumption in the
  fractional space $H^s(\Torus)$. For this reason we split the proof in two parts but in
  all cases the main step is a proper estimate of the integral
\begin{equation*}
  \mathcal{I}_{\varepsilon}:=    \int_0^t\int_{\Torus}|u_\varepsilon\otimes u_\varepsilon-(u\otimes
    u)_\varepsilon| |\nabla u_\varepsilon|\,dx d\tau.
\end{equation*}
\medskip
%
%
\textbf{The case $\frac{5}{6}\le s <1$}. Applying H\"older inequality and the convex
interpolation inequality in Lebesgue spaces we get
\begin{equation*}
\begin{split}
\mathcal{I}_{\varepsilon}&  \le \int_0^t\|u_\varepsilon\otimes u_\varepsilon-(u\otimes
  u)_\varepsilon\|_{2} \,\|\nabla u_\varepsilon\|_{2}\,d\tau
  \\
  &\le \int_0^t\|u_\varepsilon\otimes u_\varepsilon-(u\otimes
  u)_\varepsilon\|^{1-\theta}_{1}\,\|u_\varepsilon\otimes u_\varepsilon-(u\otimes
  u)_\varepsilon\|^\theta_{p}\, \|\nabla u_\varepsilon\|_{2}\,d\tau,
  \end{split}
\end{equation*}
where
\begin{equation*}
  \theta = \frac{p}{2(p-1)}.
\end{equation*}
and clearly $p\geq2$.

Note that both the $L^1$-norm of $u_{\varepsilon}(\tau)\otimes u_{\varepsilon}(\tau)$ and
$(u(\tau)\otimes u(\tau))_{\varepsilon}$ can be easily estimated by using the properties of
mollifiers as follows 
\begin{equation*}
  \|u_\varepsilon\otimes u_\varepsilon-(u\otimes u)_\varepsilon\|_{1}\le
  \|u_\varepsilon\|^2_{2}
  +\|(u\otimes u)_\varepsilon\|_{1}\le c_1\|u_\varepsilon\|^2_{2}\le c_2\|u\|^2_{2}\le C,
  \end{equation*}
  hence proving an uniform bounded since $u$ is a weak solution to~\eqref{eq:Euler}.

  Next, we fix $p=\frac{5-2s}{5-4s}$, hence $\theta=\frac{5-2s}{4s}$, and we get
\begin{equation*}
\begin{split}
\mathcal{I}_{\varepsilon}  \le C\int_0^t\|u_\varepsilon\otimes u_\varepsilon-(u\otimes
  u)_\varepsilon\|^{\frac{5-2s}{4s}}_{\frac{5-2s}{5-4s}} \|\nabla u_\varepsilon\|_{2}\,d\tau,
  \end{split}
\end{equation*}

It only remains to estimate the term involving the
$L^{\frac{5-2s}{5-4s}}$-norm. Using~\eqref{diffestimate}, and the assumption
$u\in H^s(\Torus)$ we have
\begin{equation*}
  \|u_\varepsilon\otimes u_\varepsilon-(u\otimes u)_\varepsilon\|_{p}
  \le C\varepsilon^{2\alpha}[u]^2_{W^{\alpha, 2\frac{5-2s}{5-4s}}}\le
  C\varepsilon^{2\alpha}\|u\|^2_{H^s}, 
\end{equation*}
for  $\alpha=\frac{2(s-1)s}{2s-5}$, fixed in such a way that
\begin{equation*}
    \frac{1}{2}-\frac{s}{3}= \frac{1}{2p}-\frac{\alpha}{3} \implies \alpha = \frac{3-3p+2p s}{2p},
\end{equation*}
since this is the value for which the (fractional) Sobolev embedding
$H^s(\Torus)=W^{s,2}(\Torus)\xhookrightarrow{}W^{\alpha, 2p}(\Torus)$, holds true.
Putting all together, since $2\alpha\theta=1-s$ we arrive to the following estimate
  \begin{equation*}
\mathcal{I}_{\varepsilon}\leq
C\int_0^t\|u\|_{H^s}^{\frac{5-2s}{2s}}\varepsilon^{1-s}\|\nabla u_\varepsilon\|_{2}\,d\tau.
\end{equation*}
Hence by H\"older inequality with exponents $x=\frac{5}{5-2s}$ and $x'=\frac{5}{2s}$ we
get
\begin{equation*}
\mathcal{I}_{\varepsilon}\leq
C\|u\|_{L^{5/2s}(H^s)}^{\frac{5-2s}{2s}}\ \varepsilon^{1-s}
\|\nabla u_\varepsilon\|_{L^{5/2s}(L^2)}. 
\end{equation*}
Finally by using the assumptions of Theorem~\ref{thm:energyconservationfracsob} and
 Lemma~\ref{lem:lemlimsup} we get
\begin{equation*}
  \limsup_{\varepsilon\to0^+}\mathcal{I}_{\varepsilon}
  \leq
C\limsup_{\varepsilon\to0^+}\varepsilon^{1-s}
\|\nabla u_\varepsilon\|_{L^{5/2s}(L^2)}=0. 
\end{equation*}
This is enough to end the proof since $u_{\varepsilon}(t)\to u(t)$ for
almost all $t\in(0,T)$.
\bigskip

\textbf{The case $1\leq s<5/2$}. In the case $s\geq1$ the proof is a
little different since now we can estimate directly the term
$\nabla u_{\varepsilon}$. Observe also that for $s>\frac{5}{2}$, then
$H^s(\Torus)\xhookrightarrow{}W^{1,\infty}(\Torus)$ and so one recovers the
Beale-Kato-Majda criterion for regularity if $u\in L^1(0,T;H^s(\Torus))$,
$s>\frac{5}{2}$. 

We first recall that, for $1\leq s<\frac{3}{2}$ have the
following embedding  $H^s(\Torus)\xhookrightarrow{}W^{1,p}(\Torus)$, where $p$ is such that 
$    \frac{1}{2}+\frac{1-s}{3}=\frac{1}{p}$,
and so
\begin{equation}\label{eq:pp'}
    p= \frac{6}{5-2s} \qquad \text{and} \qquad p'= \frac{p}{p-1}= \frac{6}{1+2s}.
\end{equation}
We distinguish two further sub-cases.  \smallskip

\textbf{The sub-case $1\le s < \frac{3}{2}$}. With this position we
have $2\le p < 3$. Moreover the bound $s<\frac{3}{2}$, gives even the
following embedding $H^{s}(\Torus)\xhookrightarrow{}L^{p^*}(\Torus)$, where
\begin{equation*}
    p^* = \frac{6}{3-2s}.
\end{equation*}
Applying to $\mathcal{I}_{\varepsilon}$ H\"older inequality with
conjugate exponents $p$ and $p'$, and an interpolation inequality with suitable exponents
(always possible since $p^*/2>p'$) we get
\begin{equation*}
\begin{split}
  \mathcal{I}_{\varepsilon}&\le \int_0^t\|u_\varepsilon\otimes u_\varepsilon-(u\otimes
  u)_\varepsilon\|^{1-\theta}_{1}\|u_\varepsilon\otimes u_\varepsilon-(u\otimes
  u)_\varepsilon\|^\theta_{\frac{p^*}{2}} \|\nabla u_\varepsilon\|_{p}\,d\tau,
\end{split}
\end{equation*}
with  $\theta$ satisfying the following equality
  $  \frac{1+2s}{6}= 1-\theta + \frac{\theta}{\frac{p^*}{2}}$, hence 
\begin{equation*}
    \theta = \frac{5-2s}{4s}.
\end{equation*}
As in the previous case we have 
$    \|u_\varepsilon\otimes u_\varepsilon-(u\otimes u)_\varepsilon\|_{1}\le C$.

%
Moreover, by $H^s(\Torus)\xhookrightarrow{}W^{1,p}(\Torus)$ and H\"older inequality with
exponents $5/(5-2s)$ and $5/2s$ we have
\begin{equation*}
  \mathcal{I}_{\varepsilon}\leq C \|u_\varepsilon\otimes u_\varepsilon-(u\otimes
  u)_\varepsilon\|_{L^{\frac{5}{4s}}(0,T;L^{\frac{p^*}{2}})}^{\frac{5-2s}{4s}} \|u\|_{L^\frac{5}{2s}(0,T;H^s)}.
\end{equation*}
Next,  we observe that 
\begin{equation*}
  u_\varepsilon\otimes u_\varepsilon-(u\otimes  u)_\varepsilon=u_\varepsilon\otimes
  (u_\varepsilon-u)
  +(u_{\varepsilon}-u)\otimes
  u+u\otimes u-(u\otimes
  u)_\varepsilon,
\end{equation*}
and, since $u\in L^{\frac{5}{2s}}(0,T;H^s(\Torus))$ implies
$u\otimes u\in L^{\frac{5}{4s}}(0,T;L^{\frac{p^*}{2}}(\Torus))$, then
estimate~\eqref{eq:sharpestimate} from Lemma~\ref{lem:lemlimsup} implies again that 
$\limsup_{\varepsilon\to0}\mathcal{I}_{\varepsilon}=0$, ending the proof.

\bigskip

\textbf{The sub-case $\frac{3}{2}\le s< \frac{5}{2}$}. Again we apply H\"older
inequality with conjugate exponents $p$ and $p'$ defined in~\eqref{eq:pp'} and interpolating
the $L^{p'}$-norm between $1$ and $q/2=3$ we get
\begin{equation*}
\begin{split}
  \mathcal{I}_{\varepsilon} &\le \int_0^t\|u_\varepsilon\otimes u_\varepsilon-(u\otimes
  u)_\varepsilon\|^{\frac{2s-1}{4}}_{1}\|u_\varepsilon\otimes u_\varepsilon-(u\otimes
  u)_\varepsilon\|^{\frac{5-2s}{4}}_{3} \|\nabla u_\varepsilon\|_{p}\,d\tau.
\end{split}
\end{equation*}
We use the same control as before for the $L^1(\Torus)$ norm and
%
considering the interpolation of $H^1(\Torus)$ between $L^2(\Torus)$ and $H^s(\Torus)$ and
the uniform $L^2$-bound, we have
\begin{equation*}
    \|u\|_{H^1}\le \|u\|^{1-1/s}_{L^2}\|u\|_{H^s}^{1/s}\le C \|u\|_{H^s}^{ 1/s}.
\end{equation*}
This implies also that from $u\in L^{\frac{5}{2s}}(0,T;H^s(\Torus))$ it follows --by
interpolation with $L^\infty(0,T;L^2(\Torus))$-- that
$u\in L^{\frac{5}{2s}}(0,T;H^1(\Torus))\hookrightarrow L^{\frac{5}{2s}}(0,T;L^6(\Torus))$,
hence $u\otimes u \in L^{\frac{5}{4s}}(0,T;L^3(\Torus))$. Hence by using the H\"older
inequality we get
\begin{equation*}
  \mathcal{I}_{\varepsilon}\leq C \|u_\varepsilon\otimes u_{\varepsilon}-(u\otimes
  u)_\varepsilon\|_{L^{\frac{5}{4s}}(L^{3})}^{\frac{5-2s}{4s}}\ \|u\|_{L^\frac{5}{2s}(H^s)},
\end{equation*}
and Lemma~\ref{lem:lemlimsup} implies again that
$\limsup_{\varepsilon\to0}\mathcal{I}_{\varepsilon}=0$. This ends the proof of the
conservation of energy.
\end{proof}

After having finished the proof of the criterion for energy conservation in fractional
spaces, we can pass to prove to a criterion for energy conservation who employs
vorticity/velocity in a sort of  ``geometric'' special situation. This should be compared
with the results in \cite{Der2020} where an ``analytic'' combination of the two quantities is
considered.

\begin{proof}[Proof of Theorem~\ref{thm:Beltramifieldthm}]
  Let $u\in L^\infty(0,T;L^2(\Torus))$ be a weak solution to the Euler
  equation~\eqref{eq:Euler} and let us consider $\lambda \in L^\beta(0,T;H^\tau(\Torus))$,
  for some $\beta\ge 1$ and $\tau \in \mathbb R$. Moreover, we are assuming that $u$ is a
  Beltrami field, i.e. its curl can be written as the product of $\lambda$ and itself as
  in~\eqref{eq:Beltrami}. We want to apply
  Proposition~\ref{prop:prodnonnegative}-\ref{prop:prodnegative} in order to infer sharp
  regularity for the vorticity $\omega$, which --in turn-- would give additional regularity
  for the velocity $u$. Doing so, possibly iterating, we try to show that $u$ belongs to
  some of the spaces as those in the hypotheses of
  Theorem~\ref{thm:energyconservationfracsob} to have conservation of energy.
  
  Following the notation of Propositions~\ref{prop:prodnonnegative}
  and~\ref{prop:prodnegative}, we have $s_1=0$, $s_2=\tau$ and, in both statements, it is
  required $s_i\ge s$, $i=1,2$, which gives
  \begin{equation*}
    s \le 0.
  \end{equation*}
  Moreover, a further requirement is that 
  \begin{equation*}
    s < \tau -\frac{3}{2}.
  \end{equation*}
  
  We have to distinguish different cases.

\medskip

\textbf{The case $0\leq\tau \le \frac{3}{2}$}. In this case, we have $s<0$ and we fall in
the hypotheses of Proposition~\ref{prop:prodnegative}. 
Consequently we get
\begin{equation*}
  \omega \in L^\beta(0,T;H^{\tau-\frac{3}{2}-\varepsilon}(\Torus)),\quad\text{for any
    arbitrarily     small }\varepsilon>0,
  \end{equation*}
where the integrability in time remains unchanged since $u$ is essentially bounded in
time. 

But again, $\omega$ is the curl of $u$, so that by elliptic regularity
\begin{equation*}
  u \in L^\beta(0,T;H^{\tau-\frac{1}{2}-\varepsilon}(\Torus)),\quad\text{for any arbitrarily
    small }\varepsilon>0.
  \end{equation*}
  We first note that this will give an improvement in the known regularity for the
  velocity of a weak solution only if $\tau>1/2$, hence from now on we consider $\tau$ in
  the restricted range $\tau\in]1/2,3/2]$.  Next, we can directly apply
  Theorem~\ref{thm:energyconservationfracsob} to prove conservation of energy if,
\begin{equation*}
    \frac{5}{6}<\tau-\frac{1}{2}\leq\frac{5}{2},
\end{equation*}
which holds if $\tau\in]4/3,3/2]$, and if in addition  
\begin{equation*}
    \beta > \frac{5}{2\tau-1}.
\end{equation*}
Within this range for both $\tau$ and $\beta$, the weak solution $u$ satisfies the
hypotheses of Theorem~\ref{thm:energyconservationfracsob}.

Let us now see what we can infer for smaller $\tau$, that is $\tau\in]1/2,4/3]$: We
iterate the same process with a bootstrap argument. We start the iteration of the result
on product in Sobolev spaces from
$u\in L^\beta(0,T;H^{\tau-\frac{1}{2}-\varepsilon}(\Torus))$,
$\lambda\in L^\beta(0,T;H^{\tau}(\Torus))$ and for this reason we define two sequences
$\{\beta_n\}$, $\{\sigma_{n}\}$, as follows
\begin{equation*}
  \beta_1=\beta,\qquad \text{and}\qquad \sigma_{1}=\tau-\frac{1}{2}.
\end{equation*}
Next, 
we define by recursion (which follows as a formal application of
Proposition~\ref{prop:prodnonnegative} in the limiting case
$\varepsilon=0$) 
\begin{equation*} 
  \beta_{n+1}:=\frac{\beta_n \beta}{\beta_n+\beta}\qquad \text{and}\qquad
  \sigma_{n+1}:=\min\left\{\sigma_{n},\tau,\sigma_{n}+\tau-\frac{3}{2}\right\}+1. 
\end{equation*}
\begin{rem}
  The real index $s$ for the space regularity of the velocity field after $n$ applications
  of the product theorem will be any number strictly less than $\sigma_{n}$. While
  $\beta_{n}$ will be the exact Lebesgue index with respect to the time variable.
\end{rem}
%
%
Note that, since $\tau<\frac{3}{2}$ and if $\sigma_{n}\leq\frac{3}{2}$ then
    \begin{equation*}
    \min\left\{\sigma_{n},\tau,\sigma_{n},\sigma_{n}+\tau-\frac{3}{2}\right\}=\sigma_{n}+\tau-\frac{3}{2}.
\end{equation*}
These relations imply, that
\begin{equation*}
  \beta_n=\frac{\beta}{n}\quad\text{ and}\quad
  \sigma_{n}=n(\tau-\frac{1}{2})\qquad\forall\,n\in\N. 
\end{equation*}
 By a now rigorous application of
Proposition~\ref{prop:prodnonnegative} this finally proves that %
%
\begin{equation*}
  u\in L^{\frac{\beta}{n}}(0,T;H^{n(\tau-\frac{1}{2})-\varepsilon}(\Torus)),\quad\text{for any
    arbitrarily     small }\varepsilon>0,
\end{equation*}
and this argument can be iterated as long as
$n(\tau-1/2)-\varepsilon<n(\tau-1/2)\leq\frac{3}{2}$.

Since we are considering the range $\tau\in]1/2,4/3]$ we have now
$\tau-1/2\leq\frac{5}{6}$. We then fix $n_0\in\N$ such that
\begin{equation*}
  n_0(\tau-\frac{1}{2})\leq\frac{5}{6}< (n_0+1)(\tau-\frac{1}{2}),
\end{equation*}
and iterate till reaching the regularity 
\begin{equation*}
  u\in L^{\frac{\beta}{n_0+1}}(0,T;H^{(n_0+1)(\tau-\frac{1}{2})-\varepsilon}(\Torus)),\quad \text{for
    any arbitrarily  
    small }\varepsilon>0,
\end{equation*}
which is a suitable class for energy conservation. In fact, since
$n_0(\tau-1/2)-\varepsilon<\frac{5}{6}$ the iteration is well-defined and moreover
$\tau<3/2$ implies $(n_0+1)(\tau-1/2)<5/6+1<5/2$. Finally we observe that if
$\beta>\frac{5}{2\tau-1}$, then
$\beta_{n_0+1}=\frac{\beta}{(n_0+1)}>\frac{5}{(n_0+1)(2\tau-1)}$, showing that the
hypotheses of Theorem~\ref{thm:energyconservationfracsob} are then satisfied.

\medskip

\textbf{The case $\tau > \frac{3}{2}$}. In this case note that $H^{\tau}(\Torus)
\subset L^\infty(\Torus)$,
hence we get immediately that $  \omega\in L^\beta(0,T;H^0(\Torus))$, which implies 
\begin{equation*}
  u\in L^\beta(0,T;H^1(\Torus)),
\end{equation*}
and if $\beta\geq\frac{5}{2}$ we are done, since it falls within the assumptions of
Theorem~\ref{thm:energyconservationfracsob}. Note that the result will not be improved
with a further iteration. At least in the case $\tau>\frac{5}{2}$ (but the other case
$\tau\in[3/2,5/2[$ is similar) one will get $\omega\in L^{\beta/2}(0,T;H^1(\Torus))$,
which gives
$u\in L^{\beta/2}(0,T;H^2(\Torus))\hookrightarrow L^{\beta/2}(0,T;C^{1/2}(\Torus))$, which
is an energy conservation class if $\beta\geq4$, see~\cite{Ber2023c}.

\end{proof}

\section*{Acknowledgments}
Both authors are members of INdAM GNAMPA and they are funded by MIUR within project
PRIN20204NT8W ``Nonlinear evolution PDEs, fluid dynamics and transport equations:
theoretical foundations and applications'' and MIUR Excellence, Department of Mathematics,
University of Pisa, CUP I57G22000700001.
\section*{Conflicts of interest and data availability statement}
The authors declare that there is no conflict of interest. Data sharing not applicable to
this article as no datasets were generated or analyzed during the current study.


\begin{thebibliography}{10}

\bibitem{Abe2022}
K.~Abe.
\newblock Existence of vortex rings in {B}eltrami flows.
\newblock {\em Comm. Math. Phys.}, 391(2):873--899, 2022.

\bibitem{BT2018}
C.~Bardos and E.S. Titi.
\newblock Onsager's conjecture for the incompressible {E}uler equations in
  bounded domains.
\newblock {\em Arch. Ration. Mech. Anal.}, 228(1):197--207, 2018.

\bibitem{BKM1984}
J.T. Beale, T.~Kato, and A.~Majda.
\newblock Remarks on the breakdown of smooth solutions for the $3$-{D} {E}uler
  equations.
\newblock {\em Comm. Math. Phys.}, 94(1):61--66, 1984.

\bibitem{BH2021}
A.~Behzadan and M.~Holst.
\newblock Multiplication in {S}obolev spaces, revisited.
\newblock {\em Ark. Mat.}, 59(2):275--306, 2021.

\bibitem{BY2019}
H.~Beir\~{a}o~da Veiga and J.~Yang.
\newblock On the energy equality for solutions to {N}ewtonian and
  non-{N}ewtonian fluids.
\newblock {\em Nonlinear Anal.}, 185:388--402, 2019.

\bibitem{Bel1873}
E.~Beltrami.
\newblock Sui principii fondamentali dell'idrodinamica razionale.
\newblock {\em Mem. dell'Accad. Scienze Bologna}, page 394, 1873.

\bibitem{Ber2023c}
L.~C. Berselli.
\newblock Energy conservation for weak solutions of incompressible fluid
  equations: the {H}\"{o}lder case and connections with {O}nsager's conjecture.
\newblock {\em J. Differential Equations}, 368:350--375, 2023.

\bibitem{BG2024}
L.~C. Berselli and S.~Georgiadis.
\newblock Three results on the energy conservation for the 3{D} {E}uler
  equations.
\newblock {\em NoDEA Nonlinear Differential Equations Appl.}, 31:33, 2024.


\bibitem{BDLSV2019}
T.~Buckmaster, C.~de~Lellis, L.~Sz\'{e}kelyhidi, Jr., and V.~Vicol.
\newblock Onsager's conjecture for admissible weak solutions.
\newblock {\em Comm. Pure Appl. Math.}, 72(2):229--274, 2019.

\bibitem{CCFS2008}
A.~Cheskidov, P.~Constantin, S.~Friedlander, and R.~Shvydkoy.
\newblock Energy conservation and {O}nsager's conjecture for the {E}uler
  equations.
\newblock {\em Nonlinearity}, 21(6):1233--1252, 2008.

\bibitem{CFS2010}
A.~Cheskidov, S.~Friedlander, and R.~Shvydkoy.
\newblock On the energy equality for weak solutions of the {3D}
  {N}avier-{S}tokes equations.
\newblock In {\em Contributions to current challenges in mathematical fluid
  mechanics}, Adv. Math. Fluid Mech., pages 171--175. Birkh\"auser, Basel,
  2010.

\bibitem{CET1994}
P.~Constantin, W.~E, and E.S. Titi.
\newblock Onsager's conjecture on the energy conservation for solutions of
  {E}uler's equation.
\newblock {\em Comm. Math. Phys.}, 165(1):207--209, 1994.

\bibitem{DLS2009}
C.~De~Lellis and Jr.~L. Sz{\'e}kelyhidi.
\newblock The {E}uler equations as a differential inclusion.
\newblock {\em Ann. of Math. (2)}, 170(3):1417--1436, 2009.

\bibitem{Der2020}
L.~De~Rosa.
\newblock On the helicity conservation for the incompressible {E}uler
  equations.
\newblock {\em Proc. Amer. Math. Soc.}, 148(7):2969--2979, 2020.

\bibitem{DR2000}
J.~Duchon and R.~Robert.
\newblock Inertial energy dissipation for weak solutions of incompressible
  {E}uler and {N}avier-{S}tokes equations.
\newblock {\em Nonlinearity}, 13(1):249--255, 2000.

\bibitem{EP2016}
A.~Enciso and D.~Peralta-Salas.
\newblock Beltrami fields with a nonconstant proportionality factor are rare.
\newblock {\em Arch. Ration. Mech. Anal.}, 220(1):243--260, 2016.

\bibitem{Eyi1994}
G.~L. Eyink.
\newblock Energy dissipation without viscosity in ideal hydrodynamics. {I}.
  {F}ourier analysis and local energy transfer.
\newblock {\em Phys. D}, 78(3-4):222--240, 1994.

\bibitem{Fri1995}
U.~Frisch.
\newblock {\em Turbulence, The {L}egacy of {A}.{N}.~{K}olmogorov}.
\newblock Cambridge University Press, Cambridge, 1995.

\bibitem{GLS2019}
N.~R. Gauger, A.~Linke, and P.~W. Schroeder.
\newblock On high-order pressure-robust space discretisations, their advantages
  for incompressible high {R}eynolds number generalised {B}eltrami flows and
  beyond.
\newblock {\em SMAI J. Comput. Math.}, 5:89--129, 2019.

\bibitem{Ise2018}
P.~Isett.
\newblock A proof of {O}nsager's conjecture.
\newblock {\em Ann. of Math. (2)}, 188(3):871--963, 2018.

\bibitem{LWY2023}
J.~Liu, Y.~Wang, and Y.~Ye.
\newblock Energy conservation of weak solutions for the incompressible {E}uler
  equations via vorticity.
\newblock {\em J. Differential Equations}, 372:254--279, 2023.

\bibitem{NgNgT2019}
Q.-H. Nguyen, P.-T. Nguyen, and B.~Q. Tang.
\newblock Energy equalities for compressible {N}avier-{S}tokes equations.
\newblock {\em Nonlinearity}, 32(11):4206--4231, 2019.

\bibitem{Ons1949}
L.~Onsager.
\newblock Statistical hydrodynamics.
\newblock {\em Nuovo Cimento (9)}, 6(Supplemento, 2 (Convegno Internazionale di
  Meccanica Statistica)):279--287, 1949.

\bibitem{Tri1978}
H.~Triebel.
\newblock {\em Interpolation theory, function spaces, differential operators},
  volume~18 of {\em North-Holland Mathematical Library}.
\newblock North-Holland Publishing Co., Amsterdam-New York, 1978.

\bibitem{Trk1919}
V.~Trkal.
\newblock A note on the hydrodynamics of viscous fluids.
\newblock {\em Czech J. Phys.}, 44(2):97--106, 1994.
\newblock English translation of \text{\v{C}asopis P\v{e}st.} Mat. 48 (1919)
  302--311.

\bibitem{WWWY2023}
Y.~Wang, W.~Wei, G.~Wu, and Y.~Ye.
\newblock On the energy and helicity conservation of the incompressible {E}uler
  equations.
\newblock Technical Report 2307.08322v1, arXiv, 2023.


\end{thebibliography}
\def\ocirc#1{\ifmmode\setbox0=\hbox{$#1$}\dimen0=\ht0 \advance\dimen0
  by1pt\rlap{\hbox to\wd0{\hss\raise\dimen0
  \hbox{\hskip.2em$\scriptscriptstyle\circ$}\hss}}#1\else {\accent"17 #1}\fi}
  \def\polhk#1{\setbox0=\hbox{#1}{\ooalign{\hidewidth
  \lower1.5ex\hbox{`}\hidewidth\crcr\unhbox0}}} \def\cprime{$'$}

\end{document}